\documentclass[12pt,a4paper]{article}
\usepackage{amsmath}
\usepackage{amssymb}
\usepackage{t1enc}
\usepackage[latin1]{inputenc}
\usepackage[english]{babel}
\usepackage{amsfonts}
\usepackage{latexsym}
\usepackage{bm}
\usepackage{setspace}
\usepackage[toc,page]{appendix}
\usepackage{graphicx}
\usepackage{enumerate}
\usepackage[numbers,sort&compress]{natbib}
\usepackage{tikz}

\usepackage{theoremref}
\usepackage{lineno}

\usepackage{mathtools}

\newtheorem{theorem}{Theorem}

\newtheorem{lemma}{Lemma}
\newtheorem{cor}{Corollary}

\newtheorem{construction}{Construction}

\allowdisplaybreaks

\makeatletter
\newcommand*{\Fbar}{}%
\DeclareRobustCommand*{\Fbar}{%
  \mathpalette\@Fbar{}%
}
\newcommand*{\@Fbar}[2]{%
  \sbox0{$#1\mathrm{F}\m@th$}%
  \sbox2{$#1F\m@th$}%
  \rlap{%
    \hbox to\wd2{%
      \hfill
      $\overline{%
        \vrule width 0pt height\ht0 %
        \kern\wd0 %
      }$%
    }%
  }%
  \copy2 %
}
\makeatother

\begin{document}
\title{Proof of a conjecture on isolation of graphs with a universal vertex
\medskip\medskip
}

\author{Peter Borg \\[2mm]
\normalsize Department of Mathematics, Faculty of Science, University of Malta, Malta \\
\normalsize \texttt{peter.borg@um.edu.mt}
\and
Alastair Farrugia \\ 
\normalsize \texttt{aafarrugia@gmail.com}
}

\date{}
\maketitle

\begin{abstract}
A copy of a graph $F$ is called an $F$-copy. For any graph $G$, the $F$-isolation number of $G$, denoted by $\iota(G,F)$, is the size of a smallest subset $D$ of the vertex set of $G$ such that the closed neighbourhood $N[D]$ of $D$ in $G$ intersects the vertex sets of the $F$-copies contained by $G$ (equivalently, $G-N[D]$ contains no $F$-copy). Thus, $\iota(G,K_1)$ is the domination number $\gamma(G)$ of $G$, and $\iota(G,K_2)$ is the vertex-edge domination number of $G$. Settling a conjecture of Zhang and Wu, the present first author proved that if $F$ is a $k$-edge graph, $\gamma(F) = 1$ (that is, $F$ has a vertex that is adjacent to all the other vertices of $F$), and $G$ is a connected $m$-edge graph, then $\iota(G,F) \leq \frac{m+1}{k+2} $ unless $G$ is an $F$-copy or $F$ is a $3$-vertex path and $G$ is a $6$-cycle. The graphs that attain the bound were determined for $0 \leq k \leq 1$ by Fenech, Kaemawichanurat and the present first author, and for $k = 2$ by Zhang and Wu. In this paper, we determine them for $k \geq 3$. Our result was also conjectured by Zhang and Wu.
\end{abstract}

\section{Introduction} \label{Introsection}
For standard terminology in graph theory, we refer the reader to \cite{West}. Most of the notation and terminology used here is defined in \cite{Borg}. The set of positive integers is denoted by $\mathbb{N}$. For $n \in \{0\} \cup \mathbb{N}$, $[n]$ denotes the set $\{i \in \mathbb{N} \colon i \leq n\}$. Note that $[0]$ is the empty set $\emptyset$. Arbitrary sets and graphs are taken to be finite. For a set $X$, ${X \choose 2}$ denotes the set of $2$-element subsets of $X$. 
Every graph $G$ is taken to be \emph{simple}, that is, its vertex set $V(G)$ and edge set $E(G)$ satisfy $E(G) \subseteq {V(G) \choose 2}$. We may represent an edge $\{v,w\}$ by $vw$. We call $G$ an \emph{$n$-vertex graph} if $|V(G)| = n$. We call $G$ an \emph{$m$-edge graph} if $|E(G)| = m$. For a vertex $v$ of $G$, $N_{G}(v)$ denotes the set of neighbours of $v$ in $G$, $N_{G}[v]$ denotes the closed neighbourhood $N_{G}(v) \cup \{ v \}$ of $v$, and $d_{G}(v)$ denotes the degree $|N_{G} (v)|$ of $v$. 
For a subset $X$ of $V(G)$, $N_G[X]$ denotes the closed neighbourhood $\bigcup_{v \in X} N_G[v]$ of $X$, $G[X]$ denotes the subgraph of $G$ induced by $X$, and $G - X$ denotes the graph obtained by deleting $X$ from $G$. Thus, $G[X] = (X, E(G) \cap {X \choose 2})$ and $G - X = G[V(G) \setminus X]$. We may abbreviate $G - \{x\}$ to $G - x$. Where no confusion arises, the subscript $G$ may be omitted from notation that uses it. 

For $n \geq 1$, the graphs $([n], {[n] \choose 2})$, $([n], \{\{1, i\} \colon i \in [n] \setminus \{1\}\})$ and $([n], \{\{i,i+1\} \colon i \in [n-1]\})$ are denoted by $K_n$, $K_{1,n-1}$ and $P_n$, respectively. For $n \geq 3$, $C_n$ denotes the graph $([n], \{\{1,2\}, \{2,3\}, \dots, \{n-1,n\}, \{n,1\}\})$. A copy of $K_n$ is called an \emph{$n$-clique} or a \emph{complete graph}, a copy of $K_{1,n}$ is called an \emph{$n$-star}, a copy of $P_n$ is called an \emph{$n$-path} or simply a path, and a copy of $C_n$ is called an \emph{$n$-cycle} or simply a cycle. 

If $D \subseteq V(G) = N[D]$, then $D$ is called a \emph{dominating set of $G$}. The size of a smallest dominating set of $G$ is called the \emph{domination number of $G$} and is denoted by $\gamma(G)$. If $H$ is a copy of a graph $F$, then we call $H$ an \emph{$F$-copy} and we write $H \simeq F$. If $\mathcal{F}$ is a set of graphs and $H \simeq F \in \mathcal{F}$, then we call $H$ an \emph{$\mathcal{F}$-graph}. If $H$ is a subgraph of $G$, then we say that \emph{$G$ contains $H$}. A subset $D$ of $V(G)$ is called an \emph{$\mathcal{F}$-isolating set of $G$} if $N[D]$ intersects the vertex sets of the $\mathcal{F}$-graphs contained by $G$. Thus, $D$ is an $\mathcal{F}$-isolating set of $G$ if and only if $G - N[D]$ contains no $\mathcal{F}$-graph. If $D$ is a smallest $\mathcal{F}$-isolating set of $G$, then $D$ is called a \emph{minimum $\mathcal{F}$-isolating set of $G$}, and $|D|$ is called the \emph{$\mathcal{F}$-isolation number of $G$} and is denoted by $\iota(G, \mathcal{F})$. 
If $\mathcal{F} = \{F\}$, then we may replace $\mathcal{F}$ in these defined terms and notation by $F$. 
Clearly, $D$ is a $K_1$-isolating set of $G$ if and only if $D$ is a dominating set of $G$. Thus, $\gamma(G) = \iota(G, K_1)$.

The study of isolating sets was introduced by Caro and Hansberg~\cite{CaHa17}. It is a natural generalization of the study of dominating sets \cite{C, CH, HHS, HHS2, HL, HL2}. Consider any connected $n$-vertex graph $G$. Ore \cite{Ore} showed that 
\begin{equation} \gamma(G) \leq \frac{n}{2} \label{Orebd}
\end{equation}
unless $G \simeq K_1$ (see \cite{HHS}). This is one of the earliest results in this field. While deleting the closed neighbourhood of a dominating set yields the graph with no vertices, deleting the closed neighbourhood of a $K_2$-isolating set yields a graph with no edges. In the literature, a $K_2$-isolating set is also called a \emph{vertex-edge dominating set}. Caro and Hansberg~\cite{CaHa17} proved that 
\begin{equation} \iota(G, K_2) \leq \frac{n}{3} \label{CHbd}
\end{equation}
unless $G \simeq K_2$ or $G \simeq C_5$. This was independently proved by \.{Z}yli\'{n}ski \cite{Z} and solved a problem in \cite{BCHH} (see also \cite{LMS, BG}). Let $\mathcal{C}$ be the set of cycles. Solving one of the problems posed by Caro and Hansberg \cite{CaHa17}, the first author~\cite{Borg} proved that 
\begin{equation} \iota(G,\mathcal{C}) \leq \frac{n}{4} \label{Borgcyclebound}
\end{equation} 
unless $G \simeq K_3$. A special case of this result is that $\iota(G,K_3) \leq \frac{n}{4}$ unless $G \simeq K_3$. Solving another problem in \cite{CaHa17}, Fenech, Kaemawichanurat and the first author~\cite{BFK} proved that 
\begin{equation} \iota(G, K_k) \leq \frac{n}{k+1} \label{BFKbound}
\end{equation} 
unless $G \simeq K_k$ or $k = 2$ and $G \simeq C_5$. Inequalities (\ref{Orebd}) and (\ref{CHbd}) are the cases $k = 1$ and $k = 2$, respectively. Each of the bounds in (\ref{Orebd})--(\ref{BFKbound}) is attainable if it is an integer. The graphs attaining the bound in (\ref{BFKbound}) are determined in \cite{FJKR, PX} for $k = 1$, in \cite{BG, LMS} for $k = 2$, in \cite{CCZ} for $k = 3$, and in \cite{CCZ2} for $k \geq 4$. The graphs attaining the bound in (\ref{Borgcyclebound}) are determined in \cite{CCZ}. In \cite{CaHa17}, it was also shown that $\iota(G,K_{1,k}) \leq \frac{n}{k+1}$. For $k \geq 1$, let $\mathcal{F}_{0,k} = \{K_{1,k}\}$, let $\mathcal{F}_{1,k}$ be the set of regular graphs of degree at least $k-1$, let $\mathcal{F}_{2,k}$ be the set of graphs whose chromatic number is at least $k$, and let $\mathcal{F}_{3,k} = \mathcal{F}_{0,k} \cup \mathcal{F}_{1,k} \cup \mathcal{F}_{2,k}$ (see \cite{Borgrsc, Borgrc2}). In \cite{Borgrsc}, the first author proved that for each $i \in \{0, 1, 2, 3\}$, 
\begin{equation} \iota(G,\mathcal{F}_{i,k}) \leq \frac{n}{k+1} \label{Borggenbound}
\end{equation}
unless $G \simeq K_k$ or $k = 2$ and $G \simeq C_5$, and that if $i \in \{1, 2, 3\}$ and the bound is an integer, then the bound is attainable. This generalizes (\ref{Orebd})--(\ref{BFKbound}) as $\mathcal{C} \subseteq \mathcal{F}_{1,3}$ and $K_k \in \mathcal{F}_{1,k} \cap \mathcal{F}_{2,k}$. It is worth mentioning that domination and isolation have been particularly investigated for maximal outerplanar graphs \cite{BK, BK2, CaWa13, CaHa17, Ch75, DoHaJo16, DoHaJo17, HeKa18, LeZuZy17, Li16, MaTa96, KaJi, To13}, mostly due to connections with Chv\'{a}tal's Art Gallery Theorem \cite{Ch75}.

Consider any connected $m$-edge graph $G$. Fenech, Kaemawichanurat and the first author~\cite{BFK2} also 
proved that, analogously to (\ref{BFKbound}), 
\begin{equation} \iota(G, K_k) \leq \frac{m+1}{{k \choose 2} + 2} \label{BFKbound2}
\end{equation} 
unless $G \simeq K_k$, and they also determined the graphs that attain the bound. Let $\mathcal{F}_{3,k}$ now be $\mathcal{F}_{1,k} \cup \mathcal{F}_{2,k}$. The first author \cite{Borgrc2} proved that, analogously to (\ref{Borggenbound}), for each $i \in \{1, 2, 3\}$, 
\begin{equation} \iota(G, \mathcal{F}_{i,k}) \leq \frac{m+1}{{k \choose 2} + 2} \label{Borggenbound2}
\end{equation} 
unless $G \simeq K_k$, and that if the bound is an integer, then it is attainable. This generalizes (\ref{BFKbound2}) and immediately yields $\iota(G,\mathcal{C}) \leq \frac{m+1}{5}$ if $G \not\simeq K_3$. 

In \cite{Borgdomv}, the first author proved Theorem~\ref{mainresult}, which generalizes (\ref{BFKbound2}) in another direction and verified a conjecture in \cite{ZW}. In order to state the result, we reproduce \cite[Construction~1]{Borgdomv}, which is a generalization of \cite[Construction 1.2]{BFK2} and a slight variation of the construction of $B_{n,F}$ in \cite{Borg}.

\begin{construction}[\cite{Borgdomv}] \label{const2} \emph{Consider any $m, k \in \{0\} \cup \mathbb{N}$ and any connected $k$-edge graph $F$, where $F \simeq K_1$ if $k = 0$ (that is, $V(F) \neq \emptyset$). By the division algorithm, there exist $q, r \in \{0\} \cup \mathbb{N}$ such that $m+1 = q(k+2) + r$ and $0 \leq r \leq k+1$. Let $Q_{m,k}$ be a set of size $q$. If $q \geq 1$, then let $v_1, \dots, v_q$ be the elements of $Q_{m,k}$, let $F_1, \dots, F_q$ be copies of $F$ such that the $q+1$ sets $V(F_1), \dots, V(F_q)$ and $Q_{m,k}$ are pairwise disjoint, and for each $i \in [q]$, let $w_i \in V(F_i)$, and let $G_i$ be the graph with $V(G_i) = \{v_i\} \cup V(F_i)$ and $E(G_i) = \{v_iw_i\} \cup E(F_i)$. If either $q = 0$, $T$ is the null graph $(\emptyset, \emptyset)$, and $G$ is a connected $m$-edge graph $T'$, or $q \geq 1$, $T$ is a tree with vertex set $Q_{m,k}$ (so $|E(T)| = q-1$), $T'$ is a connected $r$-edge graph with $V(T') \cap \bigcup_{i=1}^q V(G_i) = \{v_q\}$, and $G$ is a graph with $V(G) = V(T') \cup \bigcup_{i=1}^q V(G_i)$ and $E(G) = E(T) \cup E(T') \cup \bigcup_{i=1}^q E(G_i)$, then we say that $G$ is an \emph{$(m,F)$-special graph} with \emph{quotient graph $T$} and \emph{remainder graph $T'$}, and for each $i \in [q]$, we call $G_i$ an \emph{$F$-constituent of $G$}, and we call $v_i$ the \emph{$F$-connection of $G_i$ in $G$}. We say that an $(m,F)$-special graph is \emph{pure} if its remainder graph has no edges (\cite[Figure~1]{BFK2} is an illustration of a pure $(71, K_5)$-special graph). Clearly, an $(m,F)$-special graph is a connected $m$-edge graph.}
\end{construction}

If $v \in V(F) = N_F[v]$, then $v$ is called a \emph{dominating vertex of $F$} or a \emph{universal vertex of $F$}. Trivially, $\gamma(F) = 1$ if and only if $F$ has a dominating vertex. Also, $F$ is connected if it has a dominating vertex. If either $G \simeq F$, or $F \simeq K_{1,2}$ ($\simeq P_3$) and $G \simeq C_6$, then we say that $(G,F)$ is \emph{special}. We can now state the first author's result in \cite{Borgdomv}.

\begin{theorem}[\cite{Borgdomv}] \label{mainresult} If $F$ is a $k$-edge graph with $\gamma(F) = 1$, $G$ is a connected $m$-edge graph, and $(G,F)$ is not special, then 
\[\iota(G, F) \leq \frac{m+1}{k+2}.\]
Moreover, equality holds if $G$ is a pure $(m,F)$-special graph.
\end{theorem}
This proved Conjecture~4.4 of the recent paper \cite{ZW}, in which Zhang and Wu treated the case where $F$ is a $k$-star. Taking $\ell = |V(F)|$, note that $\ell - 1 \leq k \leq {\ell \choose 2}$. Clearly, $F \simeq K_{1,k}$ if $k = \ell - 1$, and $F \simeq K_{\ell}$ if $k = {\ell \choose 2}$. Thus, the two extremes of the bound in Theorem~\ref{mainresult} are the bound in \cite{ZW} and the bound in (\ref{BFKbound2}).

For a graph $F$, the \emph{complement of $F$} is the graph $(V(F), {V(F) \choose 2} \setminus E(F))$ and is denoted by $\Fbar$. For $e \in {V(F) \choose 2}$, $F+e$ denotes the graph $(V(F), E(F) \cup \{e\})$.

The graphs attaining the bound in Theorem~\ref{mainresult} are determined in \cite{BFK2} for $0 \leq k \leq 1$ (as then $F \simeq K_{k+1}$), and in \cite{ZW} for $k = 2$ (as then $F \simeq K_{1,2}$). In this paper, we show that for $k \geq 3$, the graphs attaining the bound are the pure $(m,F)$-special graphs and the $\{F+e \colon e \in E(\Fbar)\}$-graphs. Thus, we have the following result, which is proved in Section~\ref{Proofsection}.

\begin{theorem} \label{extremal}
Let $k \geq 3$, and let $F$ be a $k$-edge graph with $\gamma(F) = 1$. A connected $m$-edge graph $G$ attains the bound in Theorem~\ref{mainresult} if and only if $G$ is a pure $(m,F)$-special graph or a copy of $F+e$ for some $e \in E(\Fbar)$.
\end{theorem}
Together with Theorem~\ref{mainresult}, this proves Conjecture~4.4 in the preprint \cite{ZWpreprint} of \cite{ZW} (Theorem~\ref{extremal} was conjectured in \cite{ZWpreprint} only).

The next two sections provide preliminary results for the proof of Theorem~\ref{extremal}. We conclude this section with additional definitions and notation that will be used in the subsequent sections.

A \emph{component of $G$} is a maximal connected subgraph of $G$. Clearly, the components of $G$ are pairwise vertex-disjoint, and their union is $G$.

If $F$ and $G$ are graphs such that $G$ contains no copy of $F$, then $G$ is said to be \emph{$F$-free}.

For $t \geq 1$, a $t$-vertex path $(\{a_1, \dots, a_t\}, \{a_ia_{i+1} \colon i \in [t-1]\})$ will be represented by $(a_1, \dots, a_t)$. 

\section{The case $m \leq 2k+3$} \label{TwoCopiesSection}

In this section, we settle the case $m \leq 2k+3$ of Theorem~\ref{extremal}. The case where $m = 2k+3$ and $\iota(G,F) = 2$ needs special treatment and is addressed by the following lemma.

\begin{lemma}\label{TwoCopies}
If $k \geq 3$, $F$ is a $k$-edge graph with $\gamma(F)=1$, $F_1$ and $F_2$ are vertex-disjoint copies of $F$, and $G$ is a connected $(2k+3)$-edge graph containing $F_1$ and $F_2$, then either $\iota(G,F) = 1$ or $G$ is a pure $(2k+3, F)$-special graph.
\end{lemma}
\textbf{Proof.} 
Let $\ell = |V(F)|$. Since $k \geq 3$, $\ell \geq 3$. For each $i \in \{1, 2\}$, $F_i$ has a dominating vertex $x_i$, so $d_{F_i}(x_i) = \ell - 1$. Let $X_1 = V(F_1)$, $X_2 = V(F_2)$, $Y_1 = E(F_1)$ and $Y_2 = E(F_2)$. Let $X = V(G) \setminus (X_1 \cup X_2)$ and $Y = E(G) \setminus (Y_1 \cup Y_2)$. Then, $|Y| = 3$. 

Let $v_F$ be a dominating vertex of $F$. We will often use the following immediate observations. Since $N_F[v_F] = V(F)$,
\begin{equation} F[N_F[v_F]] = F. \label{F_obs}
\end{equation}
For each $v \in V(F)$, we have
\begin{equation} |N_{F}[v]| \geq 2 \quad \mbox{and} \quad |E(F-N_F[v])| \leq k - d_F(v_F) = k - \ell + 1,  \label{F_obs1}
\end{equation}
because if $v$ is not a dominating vertex of $F$, then $v$ and $v_F$ are two distinct members of $N_F[v]$.

Suppose $|N(x) \cap X_i| \geq 2$ for some $i \in \{1, 2\}$ and $x \in V(G) \setminus X_i$. We may assume that $i = 2$. We have $xv_1, xv_2 \in Y$ for some $v_1, v_2 \in X_2$ with $v_1 \neq v_2$. Let $e$ be the unique member of $Y \setminus \{xv_1, xv_2\}$. Since $G$ is connected, $e$ has a vertex $u_e$ in $\{x\} \cup X_1 \cup X_2$. Let $v_e$ be the other vertex in $e$. Then, $V(G) = X_1 \cup X_2 \cup \{x, v_e\}$. Suppose $x \notin X_1$. If $e = ux$ for some $u \in X_1$, then $\{x\}$ is an $F$-isolating set of $G$. Suppose $e \neq ux$ for each $u \in X_1$. Since $G$ is connected,  $e = uw$ for some $u \in X_1$ and $w \in X_2$. Let $G_w = G-N[w]$. We have $N_{G_w}[v] = N_{F_1}[v] \setminus \{u\}$ for any $v \in X_1 \cap V(G_w)$, $|N_{G_w}[v]| \leq |(N_{F_2}[v] \setminus N_{F_2}[w]) \cup \{x\}| \leq \ell - 1$ for any $v \in X_2 \cap V(G_w)$ (by (\ref{F_obs1})), and $N(x) = \{v_1, v_2\}$. If $x \in V(G_w)$ and $N_{G_w}(x) = \{v_1, v_2\}$, then $v_1, v_2 \in X_2 \setminus N_{F_2}[w]$, so we have $\ell = |X_2| \geq |N_{F_2}[w]| + |\{v_1, v_2\}| \geq 4$ by (\ref{F_obs1}). Since $|N_{G_w}[v]| \leq \ell - 1$ for each $v \in V(G_w)$, $G_w$ contains no $F$-graph, so $\{w\}$ is an $F$-isolating set of $G$. Now suppose $x \in X_1$. Then, $V(G) = X_1 \cup X_2 \cup \{v_e\}$. Let $G_x = G-N[x]$. We have $|N_{G_x}[v]| \leq |(N_{F_1}[v] \setminus N_{F_1}[x]) \cup \{v_e\}| \leq \ell - 1$ for any $v \in X_1 \cap V(G_x)$, $|N_{G_x}[v]| \leq |(N_{F_2}[v] \setminus \{v_1, v_2\}) \cup \{v_e\}| \leq \ell - 1$ for any $v \in X_2 \cap V(G_x)$, and if $v_e \in V(G_x) \setminus (X_1 \cup X_2)$, then $|N_{G_x}[v_e]| \leq |\{u_e, v_e\}| \leq \ell - 1$. Thus, $\{x\}$ is an $F$-isolating set of $G$.

Now suppose that for each $i \in \{1, 2\}$,
\begin{equation} |N(x) \cap X_i| \leq 1 \mbox{ for each } x \in V(G) \setminus X_i. \label{fact1settled}
\end{equation}

Since $G$ is connected, $G$ contains a path with only one vertex in $X_1$ and only one vertex in $X_2$. Let $P$ be a longest path $(a_1, \dots, a_t)$ contained by $G$ such that $a_1 \in X_1$, $a_t \in X_2$ and $a_i \in X$ for each $i \in [t] \setminus \{1, t\}$. Then, $t \geq 2$. Since $E(P) \subseteq Y$ and $|Y| = 3$, $t \leq 4$. If $t = 4$, then $G$ is a pure $(2k+3, F)$-special graph. Suppose $t \leq 3$.\medskip
\\
\textbf{Case 1:} $t = 2$. Let $A_{1,2} = \{uw \in Y \colon u \in X_1, \, w \in X_2\}$. Thus, $a_1a_2 \in A_{1,2}$. Let $e_1$ and $e_2$ be the two members of $Y \setminus \{a_1a_2\}$.\medskip
\\
\textbf{Case 1.1:} \emph{$e_1 \in A_{1,2}$ or $e_2 \in A_{1,2}$.} We may assume that $e_1 \in A_{1,2}$. By (\ref{fact1settled}), $e_1 = b_1b_2$ for some $b_1 \in X_1 \setminus \{a_1\}$ and $b_2 \in X_2 \setminus \{a_2\}$. Since $G$ is connected, $e_2 \cap X_1 \neq \emptyset$ or $e_2 \cap X_2 \neq \emptyset$. We may assume that $e_2 = c_1c_2$ for some $c_1 \in X_1$. If $c_2 \in X_2$, then by (\ref{fact1settled}), we have $c_1 \notin \{a_1, b_1\}$ and $c_2 \notin \{a_2, b_2\}$. 
We may assume that $d_{F_2}(a_2) \geq d_{F_2}(b_2)$ and that if $c_2 \in X_2$, then $d_{F_2}(a_2) \geq d_{F_2}(c_2)$. Let $G' = G-N[a_1]$. If $F \simeq K_3$, then $|X_2 \setminus \{a_2\}| = 2$, $V(G') \subseteq (X_2 \setminus \{a_2\}) \cup \{c_2\}$, 
and hence, since $b_2 \in X_2 \setminus \{a_2\}$ and either $c_2 \in X_2$ (so $|V(G')| \leq 2$), or $c_2 \notin X_2$ and $b_2c_2 \notin Y_1 \cup Y_2 \cup Y = E(G)$ (so $|E(G')| \leq |{V(G') \choose 2} \setminus \{b_2c_2\}| \leq 2$), $\{a_1\}$ is an $F$-isolating set of $G$. Suppose $F \not\simeq K_3$. Since $k \geq 3$, $\ell \geq 4$. For some $a_2' \in X_2 \setminus \{a_2\}$, $a_2a_2' \in E(F_2)$. We have $|E(G'[N_{G'}[v]])| \leq |E(F_1 - N_{F_1}[a_1]) \cup \{e_1, e_2\}| \leq k - \ell + 3 \leq k-1$ for any $v \in X_1 \cap V(G')$ (by (\ref{F_obs1})), $|N_{G'}[v]| \leq |N_{F_2}[v] \setminus \{a_2\}| \leq \ell - 1$ for any $v \in (X_2 \cap V(G')) \setminus \{b_2, c_2\}$, and $|E(G'[N_{G'}[b_2]])| \leq |(E(F_2) \setminus \{a_2a_2'\}) \cup \{b_1b_2\}| \leq k$. If $c_2 \in X_2$, then $|E(G'[N_{G'}[c_2]])| \leq |(E(F_2) \setminus \{a_2a_2'\}) \cup \{c_1c_2\}| \leq k$. If $c_2 \in V(G') \setminus (X_1 \cup X_2)$, then $N_{G'}[c_2] \subseteq \{c_1, c_2\}$. Suppose that $G'$ contains an $F$-copy $F'$. Then, by (\ref{F_obs}) and the above, $F' = G'[N_{G'}[b_2]]$ or $c_2 \in X_2$ and $F' = G'[N_{G'}[c_2]]$. We may assume that $F' = G'[N_{G'}[b_2]]$. Thus, $E(G'[N_{G'}[b_2]]) = (E(F_2) \setminus \{a_2a_2'\}) \cup \{b_1b_2\}$, and hence $N_{F_2}(a_2) = \{a_2'\}$. Recall that $d_{F_2}(a_2) \geq d_{F_2}(b_2)$. We have $\ell = |V(F')| = |N_{G'}[b_2]| \leq |N_{F_2}(b_2) \cup \{b_1, b_2\}| \leq d_{F_2}(a_2) + 2 \leq 3$, contradicting $\ell \geq 4$. Thus, $G'$ contains no $F$-copy, and hence $\{a_1\}$ is an $F$-isolating set of $G$.\medskip
\\
\textbf{Case 1.2:} $e_1, e_2 \notin A_{1,2}$. Since $G$ is connected, we may assume that $e_1 = b_1b_2$ for some $b_1 \in X_1$. Thus, $b_2 \notin X_2$. Since $G$ is connected, one of (a), (b) and (c) below holds:\medskip
\\
(a) $e_2 \cap X_1 \neq \emptyset$ and $e_2 \cap X_2 = \emptyset$ (as $e_2 \notin A_{1,2}$),\medskip
\\
(b) $e_2 \cap X_2 \neq \emptyset$, $e_2 \cap X_1 = \emptyset$ (as $e_2 \notin A_{1,2}$) and $b_2 \notin e_2$ (otherwise, $G$ contains $(b_1, b_2, b_3)$ with $\{b_3\} = e_2 \cap X_2$, but this contradicts the choice of $P$ and $t=2$),\medskip
\\
(c) $e_2 \cap X_1 = \emptyset = e_2 \cap X_2$ and $b_2 \in e_2$.\medskip
\\
Clearly, if $F \simeq K_3$, then $\{a_1\}$ is an $F$-isolating set of $G$. Suppose $F \not\simeq K_3$. Since $k \geq 3$, $\ell \geq 4$. If (a) holds or (c) holds, then by the argument in Case~1.1, $\{a_1\}$ is an $F$-isolating set of $G$. Suppose that (b) holds. First, suppose $d_{F_1}(a_1) \geq 2$ or $d_{F_2}(a_2) \geq 2$. By symmetry, we may assume that $d_{F_2}(a_2) \geq 2$. By the argument in Case~1.1, $\{a_1\}$ is an $F$-isolating set of $G$. Now suppose $d_{F_1}(a_1) < 2$ and $d_{F_2}(a_2) < 2$. Then, $N_{F_1}(a_1) = \{x_1\}$ and $N_{F_2}(a_2) = \{x_2\}$. Suppose $\iota(G,F) > 1$. 
Then, taking $G_i = G-N[a_i]$ for each $i \in \{1, 2\}$, $G_i$ contains an $F$-graph $H_i$ with $V(H_i) = (X_{3-i} \setminus \{a_{3-i}\}) \cup (e_{3-i} \setminus X_{3-i})$, and by the argument in Case~1.1, $H_i = G_i[V(H_i)]$. Thus, $E(G) = E(H_1) \cup E(H_2) \cup \{x_1a_1, a_1a_2, a_2x_2\}$, meaning that $G$ is a pure $(2k+3, F)$-special graph.\medskip 
\\
\textbf{Case 2:} $t = 3$. We have $Y = \{a_1a_2, a_2a_3, e\}$ for some $e \in E(G) \setminus \{a_1a_2, a_2a_3\}$. If $e \cap \{a_1, a_2, a_3\} \neq \emptyset$, then clearly $\{a_2\}$ is an $F$-isolating set of $G$. Suppose $e \cap \{a_1, a_2, a_3\} = \emptyset$. Then, since $G$ is connected, $e \cap (X_1 \setminus \{a_1\}) \neq \emptyset$ or $e \cap (X_2 \setminus \{a_3\}) \neq \emptyset$. We may assume that $e$ contains a vertex $u$ in $X_1 \setminus \{a_1\}$. Let $w$ be the other vertex in $e$, and let $X_1' = (X_1 \setminus \{a_1\}) \cup \{w\}$ and $X_2' = (X_2 \setminus \{a_3\}) \cup \{u\}$. Let $G' = G-N[a_2]$. If $G'$ contains no $F$-copy, then $\{a_2\}$ is an $F$-isolating set of $G$. Suppose that $G'$ contains an $F$-copy $F'$. Then, $w \notin X_1$, and one of (a) and (b) below holds:\medskip
\\
(a) $F' = G[X_1']$, $u$ is the only dominating vertex of $F'$ (as $N_{F'}(w) = \{u\}$), and $N_{F_1}(a_1) = \{x_1\}$ (as $|E(F')| = |E(F_1) \cup \{e\}| - d_{F_1}(a_1)$),\medskip
\\
(b) $w \in X_2$, $F' = G[X_2']$, $w$ is the only dominating vertex of $F'$, and $N_{F_2}(a_3) = \{x_2\}$.\medskip
\\
If (a) holds and $w \notin X_2$, then $G$ is a pure $(2k+3, F)$-special graph. Suppose that (a) holds and $w \in X_2$. Let $G^* = G-N[w]$. We have $|N_{G^*}[v]| \leq |(N_{F_2}[v] \setminus N_{F_2}[w]) \cup \{a_2\}| \leq \ell - 1$ for any $v \in X_2 \cap V(G^*)$ (by (\ref{F_obs1})), $|N_{G^*}[v]| \leq |N_{F_1}[v] \setminus \{u\}| \leq \ell - 1$ for any $v \in (X_1 \setminus \{a_1\}) \cap V(G^*)$, $G^*[N_{G^*}[a_1]] = (\{x_1, a_1, a_2\}, \{x_1a_1, a_1a_2\})$ and $G^*[N_{G^*}[a_2]] = (\{a_1, a_2, a_3\}, \{a_1a_2, a_2a_3\})$. By (\ref{F_obs}) and $k \geq 3$, $G^*$ contains no $F$-copy, so $\{w\}$ is an $F$-isolating set of $G$. Now suppose that (b) holds. Then, by symmetry, we can apply an argument similar to that for the case where (a) holds and $w \in X_2$.\hfill{$\Box$}

\begin{cor}
 \label{m<=2k+3} 
Theorem~\ref{extremal} holds if $m \leq 2k+3$.
\end{cor}
\textbf{Proof.} Suppose $\iota(G,F) = \frac{m+1}{k+2}$. Then, $0 < \iota(G,F) \leq \frac{2k+4}{k+2} = 2$. Suppose $\iota(G,F) = 1$. Then, $m = k + 1$, $G$ contains an $F$-copy $H$, and $E(G) = E(H) \cup \{vw\}$ for some $v, w \in V(G)$ with $v \in V(H)$ (as $G$ is connected) and $vw \notin E(H)$. If $w \not\in V(H)$, then $G$ is a pure $(m,F)$-special graph. If $w \in V(H)$, then $G = H + vw$. Now suppose $\iota(G,F) = 2$. Then, $m = 2k+3$. Let $F_1$ be an $F$-copy contained by $G$. Let $x_1$ be a dominating vertex of $F_1$. Since $\iota(G,F) = 2$, $G-N[x_1]$ contains an $F$-copy $F_2$. Since $V(F_1) \subseteq N[x_1]$, $V(F_1) \cap V(F_2) = \emptyset$. By Lemma~\ref{TwoCopies}, $G$ is a pure $(m,F)$-special graph.~\hfill{$\Box$}

\section{Other lemmas} \label{SupportingLemmas}

In this section, we provide additional tools for the proof of Theorem~\ref{extremal}. We start with two lemmas from \cite{Borg}.

\begin{lemma}[\cite{Borg}] \label{lemma}
If $G$ is a graph, $\mathcal{F}$ is a set of graphs, $X \subseteq V(G)$, and $Y \subseteq N[X]$, then $\iota(G, \mathcal{F}) \leq |X| + \iota(G-Y, \mathcal{F})$. 
\end{lemma}
\textbf{Proof.} Let $D$ be an $\mathcal{F}$-isolating set of $G-Y$ of size $\iota(G-Y, \mathcal{F})$. Clearly, $V(F) \cap Y \neq \emptyset$ for each $\mathcal{F}$-graph $F$ that is a subgraph of $G$ and not a subgraph of $G-Y$. Since $Y \subseteq N[X]$, $X \cup D$ is an $\mathcal{F}$-isolating set of $G$. The result follows.~\hfill{$\Box$}

\begin{lemma}[\cite{Borg, Borgrsc}] \label{lemmacomp}
If $G_1, \dots, G_r$ are the distinct components of a graph $G$, and $\mathcal{F}$ is a set of connected graphs, 
then $\iota(G,\mathcal{F}) = \sum_{i=1}^r \iota(G_i,\mathcal{F})$.
\end{lemma}
\textbf{Proof.} For each $i \in [r]$, let $D_i$ be a smallest $\mathcal{F}$-isolating set of $G_i$. Consider any $\mathcal{F}$-graph $F$ contained by $G$. Since $F$ is connected, there exists some $j \in [r]$ such that $G_j$ contains $F$, so $N[D_j] \cap V(F) \neq \emptyset$. Thus, $\bigcup_{i = 1}^r D_i$ is an $\mathcal{F}$-isolating set of $G$, and hence $\iota(G, \mathcal{F}) \leq \sum_{i = 1}^r |D_i| = \sum_{i = 1}^r \iota(G_i, \mathcal{F})$. Let $D$ be a smallest $\mathcal{F}$-isolating set of $G$. For each $i \in [r]$, $D \cap V(G_i)$ is an $\mathcal{F}$-isolating set of $G_i$. We have $\sum_{i = 1}^r \iota(G_i, \mathcal{F}) \leq \sum_{i = 1}^r |D \cap V(G_i)| = |D| = \iota(G, \mathcal{F})$. The result follows.~\hfill{$\Box$}

\begin{lemma}\label{MinIsolatingSets}
Let $F$ be a connected graph with $V(F) \neq \emptyset$. Let $G$ be a pure $(|E(G)|,F)$-special graph or a copy of $F+e$ for some $e \in E(\Fbar)$. Then, 
\begin{equation} \iota(G,F) = \frac{|E(G)|+1}{|E(F)|+2}, \label{extgrapheqn}
\end{equation}
and if $G$ is a pure $(|E(G)|,F)$-special graph, then
\begin{equation} \iota(G,F) = \frac{|V(G)|}{|V(F)| + 1}. \label{extgrapheqn2}
\end{equation}
Moreover, for each $x \in V(G)$, 
\begin{enumerate}[(a)]
\item $G$ has a minimum $F$-isolating set $D$ with $x \in D$, and
\item if $G$ is a pure $(|E(G)|,F)$-special graph and $x$ is not a vertex of a quotient graph of $G$, then $\iota(G-x,F) = \iota(G,F) - 1$.
\end{enumerate}
\end{lemma}
\textbf{Proof.} The result is trivial if $G \simeq F+e$ for some $e \in E(\Fbar)$. Let $k = |E(F)|$, $\ell = |V(F)|$, $m = |E(G)|$ and $n = |V(G)|$. Suppose that $G$ is a pure $(m,F)$-special graph as in Construction~\ref{const2}. Then, $m+1 = q(k+2)$, $n = q(\ell + 1)$, and $\{v_1, \dots, v_q\}$ is an $F$-isolating set of $G$, so $\iota(G, F) \leq q = \frac{m+1}{k+2} = \frac{n}{l+1}$. If $D$ is an $F$-isolating set of $G$, then, since $G_1 - v_1, \dots, G_q - v_q$ are copies of $F$, $D \cap V(G_i) \neq \emptyset$ for each $i \in [q]$. Thus, $\iota(G, F) = q$, and hence (\ref{extgrapheqn}) and (\ref{extgrapheqn2}) are proved. 

Let $x \in V(G)$. We may assume that $x \in V(G_1)$. Since $F$ is connected, $G_1$ is connected, so $|N_{G_1}[x]| \geq 2$, and hence $|V(G_1 - N_{G_1}[x])| < |V(F_1)|$. Let $D = \{v_i \colon i \in [q] \setminus \{1\}\}$. Then, $\{x\} \cup D$ is an $F$-isolating set of $G$ of size $q$, and hence (a) is proved.

Suppose $x \neq v_1$. Then, $x \in V(F_1)$. If $q \geq 2$, then, since $T$ is connected, we have $v_1 \in N[D]$, so $D$ is an $F$-isolating set of $G-x$. Suppose $q = 1$ and $\iota(G-x) = 1$. Then, $G = G_1$, $G-x \simeq F$ and $E(G) = E(G-x) \cup \{xy\}$ for some $y \in (V(F_1) \setminus \{x\}) \cup \{v_1\}$. Thus, $(\{x\}, \emptyset)$ is a quotient graph of $G$, and hence (b) is proved.~\hfill{$\Box$}

\begin{lemma}\label{MinIsolatingSets2vertices}
If $k \geq 3$, $F$ is a $k$-edge graph with $\gamma(F)=1$, $H$ and $I$ are distinct $F$-constituents of a pure $(m,F)$-special graph $G$, $x \in V(H)$ and $y \in V(I)$, then $G$ has a minimum $F$-isolating set $D$ with $x, y \in D$.
\end{lemma}
\textbf{Proof.} We may assume that $G$ is as in Construction~\ref{const2}, $H = G_1$ and $I = G_2$. Let $D = \{x, y\} \cup \{v_i \colon i \in [q] \setminus [2]\}$. If $y \in \{v_2, w_2\}$, then $G-N[D]$ is a subgraph of $(G_1 - N_{G_1}[x]) \cup \bigcup_{i=2}^q (F_i - w_i)$, which is $F$-free as $|N_{G_1}[x]| \geq 2$. Similarly, $G-N[D]$ is $F$-free if $x \in \{v_1, w_1\}$. Suppose $x \not\in \{v_1, w_1\}$ and $y \not\in \{v_2, w_2\}$.
Let $G^* = G-N[D]$. 
Suppose that $G^*$ contains a copy $F'$ of $F$ with dominating vertex $u$. Then, $|E(G^*[N_{G^*}[u]])| \geq |E(G^*[V(F')])| \geq k$. Since $k \geq 3$ and $E(G^*[N_{G^*}[v_1]]) \subseteq \{v_1w_1, v_1v_2\}$, $u \neq v_1$. Similarly, $u \neq v_2$. Thus, $F'$ is a subgraph of $G_1 - N[x]$ or of $G_2 - N[y]$. However, $|V(G_1 - N[x])| \leq |V(F)| - 1$ and $|V(G_2 - N[y])| \leq |V(F)| - 1$, a contradiction.~\hfill{$\Box$}

\begin{lemma}\label{2PureSpecial}
If $k \geq 3$, $F$ is a $k$-edge graph with $\gamma(F)=1$, $\Gamma_1$ is a pure $(m_1,F)$-special graph, $\Gamma_2$ is a pure $(m_2,F)$-special graph, $V(\Gamma_1) \cap V(\Gamma_2) = \emptyset$, $\Gamma = (V(\Gamma_1) \cup V(\Gamma_2), E(\Gamma_1) \cup E(\Gamma_2) \cup \{x_1 x_2\})$ for some $x_1 \in V(\Gamma_1)$ and $x_2 \in V(\Gamma_2)$, $m = |E(\Gamma)|$ and $\iota(\Gamma,F) = \frac{m+1}{k+2}$, then $\Gamma$ is a pure $(m,F)$-special graph.
\end{lemma}
\textbf{Proof.} Let $\Gamma_1$ have quotient graph $T_1$, $F$-connections $v_1, \ldots, v_{q_1}$, and $F$-constituents $G_1, \ldots, G_{q_1}$, where $v_i \in V(G_i)$ for each $i \in [q_1]$. For each $i \in [q_1]$, $V(G_i) = \{v_i\} \cup V(F_i)$ for some $F$-copy $F_i$ with $v_i \notin V(F_i)$, and $v_iw_i \in E(G_i)$ for some $w_i \in V(F_i)$. Let $\Gamma_2$ have quotient graph $T_2$, and let $q_2 = |V(T_2)|$. Thus, $m_i + 1 = q_i(k+2)$ for each $i \in [2]$. Let $q = q_1 + q_2$. We have $m + 1 = m_1 + m_2 + 2 = q(k+2)$. Suppose $\iota(\Gamma,F) = \frac{m+1}{k+2}$. Then, $\iota(\Gamma,F) = q$.  





If $x_1 \in V(T_1)$ and $x_2 \in V(T_2)$, then $\Gamma$ is a pure $(m,F)$-special graph with quotient graph $(V(T_1) \cup V(T_2), E(T_1) \cup E(T_2) \cup \{x_1x_2\})$. 

Now suppose $x_1 \notin V(T_1)$ or $x_2 \notin V(T_2)$. We may assume that $x_1 \in V(F_1)$. By Lemma~\ref{MinIsolatingSets}, for each $i \in [2]$, $\Gamma_i$ has an $F$-isolating set $D_i$ such that $|D_i| = q_i$ and $x_i \in D_i$. If $x_1 = w_1$, then $\{w_i \colon i \in [q_1] \setminus \{1\}\} \cup D_2$ is an $F$-isolating set of $\Gamma$ of size $q - 1$, but this contradicts $\iota(\Gamma,F) = q$. Thus, $x_1 \in V(F_1) \setminus \{w_1\}$. Let $D_1' = \{v_i \colon i \in [q_1] \setminus \{1\}\}$ and $D = D_1' \cup D_2$. If $q_1 > 1$, then $v_1 \in N[D_1']$ (as in the proof of Lemma~\ref{MinIsolatingSets}), so $D$ is an $F$-isolating set of $\Gamma$, but this contradicts $\iota(\Gamma,F) = q$. Thus, $\Gamma_1 = G_1$ and $\Gamma-N[D] = (G_1 - x_1) \cup (\Gamma_2 - N_{\Gamma_2}[D_2])$. Since $\Gamma_2 - N_{\Gamma_2}[D_2]$ is $F$-free and $\iota(\Gamma, F) = q = 1 + |D_2|$, $\iota(G_1 - x_1, F) \neq 0 = \iota(G_1, F) - 1$. By Lemma~\ref{MinIsolatingSets} (b), $(\{x_1\}, \emptyset)$ is a quotient graph of $G_1$. Thus, $\Gamma$ is a pure $(m,F)$-special graph with quotient graph $(\{x_1\} \cup V(T_2), \{x_1 x_2\} \cup E(T_2))$.~\hfill{$\Box$}


\begin{lemma}\label{2PureSpecialMulti} Let $k \geq 3$ and let $F$ be a $k$-edge graph with $\gamma(F)=1$. Let $s \geq 2$ and let $\Gamma_1, \ldots, \Gamma_s$ be vertex-disjoint graphs such that for each $i \in [s]$, $\Gamma_i$ is a pure $(|E(\Gamma_i)|,F)$-special graph or a copy of $F+e$ for some $e \in E(\Fbar)$. For each $i \in [s] \setminus \{1\}$, let $x_i \in V(\Gamma_1)$ and $y_i \in V(\Gamma_i)$. Let $\Gamma$ be the graph with $V(\Gamma) = \bigcup_{i=1}^s V(\Gamma_i)$ and $E(\Gamma) = E(\Gamma_1) \cup \bigcup_{i=2}^{s} (E(\Gamma_i) \cup \{x_i y_i\})$. Let $m = |E(\Gamma)|$. If $\iota(\Gamma,F) = \frac{m+1}{k+2}$, then $\Gamma$ is a pure $(m,F)$-special graph.
\end{lemma}
\textbf{Proof.} We use induction on $s$. Suppose $\iota(\Gamma,F) = \frac{m+1}{k+2}$. For each $i \in [s]$, let $m_i = |E(\Gamma_i)|$. By Lemma~\ref{MinIsolatingSets}, $\iota(\Gamma_i,F) = \frac{m_i+1}{k+2}$, and if $i \neq 1$, then $\Gamma_1$ has a minimum $F$-isolating set $D_{1,i}$ containing $x_i$, and $\Gamma_i$ has a minimum $F$-isolating set $D_i$ containing $y_i$. Let $D_1 = D_{1,2}$ and $D = \bigcup_{i=1}^s D_i$. Then, $D$ is an $F$-isolating set of $\Gamma$, and
\begin{equation}
|D| = \sum_{i=1}^s \frac{m_i+1}{k+2} = \frac{m+1}{k+2} = \iota(\Gamma,F). \nonumber    
\end{equation}

Suppose that for some $j \in [s]$, $\Gamma_{j} \simeq F+e$ for some $e \in E(\Fbar)$. Suppose $j = 1$. Let $D' = D \setminus D_1$. Since $y_2 \in D_2 \subseteq D'$ and $x_2 \in N(y_2)$, $D'$ is an $F$-isolating set of $\Gamma$. We have $|D'| < |D|$, which contradicts $|D| = \iota(\Gamma,F)$. Thus, $j > 1$. We have that $(D \setminus (D_1 \cup D_j)) \cup D_{1,j}$ is an $F$-isolating set of $\Gamma$ that is smaller than $|D|$, and this again contradicts $|D| = \iota(\Gamma,F)$. 

Therefore, for each $i \in [s]$, $\Gamma_i$ is a pure $(m_i,F)$-special graph. If $s=2$, then the result is given by Lemma~\ref{2PureSpecial}. Suppose $s \geq 3$. Let $\Gamma^* = \Gamma - V(\Gamma_s)$. Let $m^* = |E(\Gamma^*)|$. Then, $m^* = m - m_s - 1$. Let $D^* = D \setminus D_s$. Since $y_j \in D_j$ for each $j \in [s] \setminus \{1\}$, $x_2 y_2, \ldots, x_{s-1} y_{s-1} \notin E(\Gamma^* - N_{\Gamma^*}[D^*])$, so $\Gamma^* - N_{\Gamma^*}[D^*]$ is a subgraph of $\bigcup_{i = 1}^{s-1} (\Gamma_i - N_{\Gamma_i}[D_i])$. Thus, $\Gamma^* - N_{\Gamma^*}[D^*]$ is $F$-free, and hence $\iota(\Gamma^*,F) \leq |D^*| = \frac{m^* + 1}{k+2}$.

Suppose that $\Gamma^*$ has an $F$-isolating set $D''$ with $|D''| < \frac{m^* + 1}{k+2}$. 
Since $y_s \in D_s$, $D'' \cup D_s$ is an $F$-isolating set of $\Gamma$. We have $\iota(\Gamma,F) \leq |D''| + |D_s| < \frac{m^* + 1}{k+2} + \frac{m_s + 1}{k+2} = \frac{m + 1}{k+2}$, which contradicts $\iota(\Gamma,F) = \frac{m + 1}{k+2}$. Thus, $\iota(\Gamma^*,F) = \frac{m^* + 1}{k+2}$. By the induction hypothesis, $\Gamma^*$ is a pure $(m^*,F)$-special graph. By applying Lemma~\ref{2PureSpecial} to $\Gamma^*$ and $\Gamma_s$, we conclude that $\Gamma$ is a pure $(m,F)$-special graph.
~\hfill{$\Box$}

\section{Proof of Theorem~\ref{extremal}} \label{Proofsection}

For a vertex $v$ of a graph $G$, let $E_G(v)$ denote the set $\{vw \colon w \in N_G(v)\}$. For $X, Y \subseteq V(G)$, let $E_G(X,Y)$ denote the set $\{xy \in E(G) \colon x \in X, \, y \in Y\}$.\\
\\
\textbf{Proof of Theorem~\ref{extremal}.} By Lemma~\ref{MinIsolatingSets}, $G$ attains the bound in Theorem~\ref{mainresult} if $G$ is a pure $(m,F)$-special graph or a copy of $F+e$ for some $e \in E(\Fbar)$. We now prove the converse.

Let $n = |V(G)|$ and $\ell = |V(F)|$. We use induction on $n$. The result is given by Corollary~\ref{m<=2k+3} if $m \leq 2k+3$. Suppose $m \geq 2k+4$. We have $\iota(G,F) = \frac{m+1}{k+2} > 2$. If $n \leq \ell$, then $\iota(G,F) \leq 1$, a contradiction. Thus, $n \geq \ell + 1$. Since $k \geq 3$, $\ell \geq 3$. Let $\mathcal{S}$ be the set of $F$-copies contained by $G$. Since $\gamma(F) = 1$, for each $S \in \mathcal{S}$, $V(S) = N_S[v_S]$ for some $v_S \in V(S)$. Let $U = \{u \in V(G) \colon V(S) = N_S[u] \mbox{ for some } S \in \mathcal{S}\}$. Let $v \in U$ such that $d_G(u) \leq d_G(v)$ for each $u \in U$. For some $F_1 \in \mathcal{S}$, $V(F_1) = N_{F_1}[v] \subseteq N_G[v]$. Thus, $d(v) \geq \ell - 1 \geq 2$. Since $\iota(G, F) > 2$, $V(G) \neq N[v]$. Let $G' = G-N[v]$ and $n' = |V(G')|$. Then, $V(G') \neq \emptyset$.

Let $\mathcal{H}$ be the set of components of $G'$. For any $H \in \mathcal{H}$ and any $x \in N(v)$ such that $xy_{x,H} \in E(G)$ for some $y_{x,H} \in V(H)$, we say that $H$ is \emph{linked to $x$} and that $x$ is \emph{linked to $H$}. Since $G$ is connected, for each $H \in \mathcal{H}$, $x_Hy_H \in E(G)$ for some $x_H \in N(v)$ and some $y_H \in V(H)$. We have 
\begin{equation} E(F_1) \subseteq E(G[N[v]]), \quad \{x_Hy_H \colon H \in \mathcal{H}\} \subseteq E_G(N(v), V(G')), \label{edgesubsets1}
\end{equation}
\begin{equation}  E(F_1) \cup \bigcup_{H \in \mathcal{H}} (E(H) \cup \{x_Hy_H\}) \subseteq E(G), \label{extremeH}
\end{equation} 
\begin{equation} m = |E(G[N[v]])| + |E_G(N(v), V(G'))| + \sum_{H \in \mathcal{H}} |E(H)| \geq k + \sum_{H \in \mathcal{H}} |E(H) \cup \{x_Hy_H\}|. \label{sizeineq_1}
\end{equation} 
Let $\mathcal{H}' = \{H \in \mathcal{H} \colon (H,F) \mbox{ is special}\}$. 
Since $k \geq 3$, we have $F \not\simeq K_{1,2}$, so
$\mathcal{H}' = \{H \in \mathcal{H} \colon H \simeq F\}$. 
By Theorem~\ref{mainresult},
\begin{equation} \iota(H, F) \leq \frac{|E(H)|+1}{k+2} \quad \mbox{for each } H \in \mathcal{H} \setminus \mathcal{H}'. \label{E(H)_eqn}
\end{equation}
\textbf{Case 1:} \emph{$\mathcal{H}' = \emptyset$.} By Lemma~\ref{lemma} (with $X = \{v\}$ and $Y = N[v]$) and Lemma~\ref{lemmacomp}, $\iota(G, F) \leq 1 + \iota(G', F) = 1 + \sum_{H \in \mathcal{H}} \iota(H, F)$. Since $\iota(G, F) = \frac{m+1}{k+2} > \sum_{H \in \mathcal{H}} \frac{|E(H)|+1}{k+2} \geq \sum_{H \in \mathcal{H}} \iota(H, F)$, we actually have 
\begin{equation} 
\iota(G, F) = 1 + \iota(G', F) = 1 + \sum_{H \in \mathcal{H}} \iota(H, F). \label{isol-components} 
\end{equation}
%
%
%
Let $\sigma = k + \sum_{H \in \mathcal{H}} (|E(H)|+1)$. 
%
%
Let
\begin{gather} A = E(G[N[v]]) \setminus E(F_1), \quad B = E_G(N(v), V(G')) \setminus \{x_H y_H: H \in \mathcal{H}\}, \nonumber \\
\mathcal{H}_1 = \left\{ H \in \mathcal{H}: \iota(H,F) = \frac{|E(H)|+1}{k+2} \right\}, \quad \mathcal{H}_2 = \left\{ H \in \mathcal{H}: \iota(H,F) < \frac{|E(H)|+1}{k+2} \right\}, \nonumber \\
a = |A|, \quad b = |B|, \quad c = \sum_{H \in \mathcal{H}} (|E(H)|+1-(k+2)\iota(H,F)). \nonumber
\end{gather}
By (\ref{E(H)_eqn}), $\mathcal{H} = \mathcal{H}_1 \cup \mathcal{H}_2$ and $c \geq 0$. We have
\begin{equation} E(G) = E(F_1) \cup A \cup B \cup \bigcup_{H \in \mathcal{H}} (E(H) \cup \{x_H y_H\}), \label{E(G)_eqn}
\end{equation}
so $m = \sigma + a + b$. We have
\begin{align} \frac{m+1}{k+2} &= \iota(G, F) = 1 + \sum_{H \in \mathcal{H}} \iota(H, F) = \frac{k+2}{k+2} - \frac{c}{k+2} + \sum_{H \in \mathcal{H}} \frac{|E(H)|+1}{k+2} \nonumber \\
&= \frac{\sigma - c +2}{k+2} = \frac{m - a - b - c + 2}{k+2}, \label{iotaineq_0}
\end{align}
so $a + b + c = 1$. Thus, one of $a$, $b$ and $c$ is $1$, and the rest are $0$.\medskip
\\
\textbf{Case 1.1} \emph{$c = 0$.} Thus, for each $H \in \mathcal{H}$, $\iota(H,F) = \frac{|E(H)|+1}{k+2}$, so by the induction hypothesis, $H$ is a pure $(|E(H)|,F)$-special graph or a copy of $F+e$ for some $e \in E(\Fbar)$, and by Lemma~\ref{MinIsolatingSets}, $H$ has a minimum $F$-isolating set $D_H$ containing $y_H$. 
Let $D = \bigcup_{H \in \mathcal{H}} D_H$ 
and $X = \{x_H \colon H \in \mathcal{H}\}$. 

Suppose $a = 1$. Then, $E(G) = E(F_1) \cup \{e^*\} \cup \bigcup_{H \in \mathcal{H}} (E(H) \cup \{x_H y_H\})$ for some $e^* \in E(G[N[v]])$, and $G - N[D] = (G[N[v]] - X) \cup \bigcup_{H \in \mathcal{H}} (H - N_H[D_H])$. By (\ref{isol-components}), $D$ is not an $F$-isolating set of $G$, so $G[N[v]] - X$ contains a copy $F_2$ of $F$. Since $E(G[N[v]]) = E(F_1) \cup \{e^*\}$, it follows that $E(G[N[v]]) = E(F_2) \cup \{vx\}$ with $\{x\} = X$ and $x \notin V(F_2)$, so $G[N[v]]$ is a pure $(k+1,F)$-special graph. By Lemma~\ref{2PureSpecialMulti}, $G$ is a pure $(m,F)$-special graph.

Now suppose $b = 1$. Then, $E(G) = E(F_1) \cup \{\tilde{e}\} \cup \bigcup_{H \in \mathcal{H}} (E(H) \cup \{x_H y_H\})$ with $\{\tilde{e}\} = B$. For some $H_0 \in \mathcal{H}$, $\tilde{e} = \tilde{x} \tilde{y}$ for some $\tilde{x} \in N(v)$ and $\tilde{y} \in V(H_0)$. Let $x_0 = x_{H_0}$, $y_0 = y_{H_0}$, $D_0 = D_{H_0}$ and $G^0 = G[N[v] \cup V(H_0)]$.

Since $\{v\} \cup D_0$ is an $F$-isolating set of $G^0$, $\iota(G^0,F) \leq 1 + \iota(H_0,F)$. Suppose $\iota(G^0,F) \leq \iota(H_0,F)$. 
Then, $G^0$ has an $F$-isolating set $D^*$ with $|D^*| \leq \iota(H_0,F)$. Since $y_H \in D_H$ for each $H \in \mathcal{H} \setminus \{H_0\}$, $D^* \cup \bigcup_{H \in \mathcal{H} \setminus \{H_0\}} D_H$ is an $F$-isolating set of $G$. This yields $\iota(G,F) \leq \sum_{H \in \mathcal{H}} \iota(H, F)$, which contradicts (\ref{isol-components}). Thus,  
\begin{equation} \iota(G^0,F) = 1 + \iota(H_0,F) = \frac{|E(F_1) \cup \{\tilde{e}\}| + 1}{k+2} + \frac{|E(H_0) \cup \{x_0y_0\}|}{k+2} = \frac{|E(G^0)|+1}{k+2}. \label{G0edges}
\end{equation} 

Suppose $\iota(H_0,F) = 1$. Then, $2 = \iota(G^0,F) = \frac{|E(G^0)|+1}{k+2}$, so $|E(G^0)| = 2k + 3$. By Lemma~\ref{TwoCopies}, $G^0$ is a pure $(2k+3, F)$-special graph, so $|V(G^0)| = 2\ell + 2$ by (\ref{extgrapheqn2}). Since $1 = \iota(H_0,F) = \frac{|E(H_0)|+1}{k+2}$, $|V(H_0)| \in \{\ell, \ell + 1\}$ by Corollary~\ref{m<=2k+3}. Since $V(G^0) = V(F_1) \cup V(H_0)$, we have $|V(G^0)| \leq 2\ell + 1$, a contradiction. 

Therefore, $\iota(H_0,F) \geq 2$, and hence $H_0$ is a pure $(|E(H_0)|,F)$-special graph. Let $T$ be the quotient graph of $H_0$, let $v_1, \ldots, v_q$ be the distinct vertices of $T$, and let $G_1, \ldots, G_q$ be the $F$-constituents of $H_0$, where $v_i \in V(G_i)$ for each $i \in [q]$. We may assume that $y_0 \in V(G_1)$. Suppose $\tilde{y} \notin V(G_1)$. Then, by Lemma~\ref{MinIsolatingSets2vertices}, $H_0$ has a minimum $F$-isolating set $D_0'$ with $y_0, \tilde{y} \in D_0'$. Thus, $D_0'$ is an $F$-isolating set of $G^0$, but this contradicts $\iota(G^0,F) = 1 + \iota(H_0,F)$. Thus, $\tilde{y} \in V(G_1)$. Let $J = G[V(F_1) \cup V(G_1)]$. Then, $E(J) = E(F_1) \cup \{x_0y_0, \tilde{x}\tilde{y}\} \cup E(G_1)$, so $J$ is connected and $|E(J)| = 2k+3$. By Theorem~\ref{mainresult}, $\iota(J,F) \leq 2$. By (\ref{extgrapheqn2}), since $|V(J)| = 2 \ell + 1$, $J$ is not a pure $(2k+3, F)$-special graph. By the induction hypothesis, $\iota(J,F) \neq 2$, so $J - N_J[v_J]$ is $F$-free for some $v_J \in V(J)$. Thus, $\{v_J, v_2, \dots, v_q\}$ is an $F$-isolating set of $G^0$ of size $q$, and hence we have $\iota(G^0,F) \leq q = \frac{|E(H_0)|+1}{k+2} = \iota(H_0, F)$, contradicting (\ref{G0edges}).
%
\medskip
\\
\textbf{Case 1.2:} \emph{$c=1$.} Thus, $E(G) = E(F_1) \cup \bigcup_{H \in \mathcal{H}} (E(H) \cup \{x_H y_H\})$, $|E(H_0)| = (k+2)\iota(H_0,F)$ for some $H_0 \in \mathcal{H}$, and $|E(H)| + 1 = (k+2)\iota(H,F)$ for each $H \in \mathcal{H} \setminus \{H_0\}$. By the induction hypothesis, for each $H \in \mathcal{H} \setminus \{H_0\}$, $H$ is a pure $(|E(H)|,F)$-special graph or a copy of $F+e$ for some $e \in E(\Fbar)$, and by Lemma~\ref{MinIsolatingSets}, $H$ has a minimum $F$-isolating set $D_H$ containing $y_H$. Let $x_0 = x_{H_0}$ and $y_0 = y_{H_0}$. Let $G^0 = G[N[v] \cup V(H_0)]$. Similarly to Case~1.1, we have
\begin{equation} \iota(G^0,F) = 1 + \iota(H_0,F) = \frac{|E(F_1) \cup \{x_0y_0\}| + 1}{k+2} + \frac{|E(H_0)|}{k+2} = \frac{|E(G^0)|+1}{k+2}. \nonumber
\end{equation} 

Suppose $\mathcal{H} \neq \{H_0\}$. Let $X = \{x_H \colon H \in \mathcal{H} \setminus \{H_0\}\}$ and $D_X = \bigcup_{H \in \mathcal{H} \setminus \{H_0\}} D_H$. Suppose $x_0 \in X$. Let $D_0$ be a minimum $F$-isolating set of $H_0$. Then, $D_0 \cup D_X$ is an $F$-isolating set of $G$ of size $\sum_{H \in \mathcal{H}} \iota(H,F)$, which contradicts (\ref{isol-components}). Thus, $x_0 \notin X$. Let $G^0_X = G^0 - X$. We have $v, x_0 \in N[v] \setminus X \subseteq V(G^0_X) = (N[v] \setminus X) \cup V(H_0)$, so $G^0_X$ is connected. Also, $\{vx \colon x \in X\} \subseteq E(G^0) \setminus E(G^0_X)$, so $|E(G^0_X)| < |E(G^0)|$. Suppose $G^0_X \not\simeq F$. By Theorem~\ref{mainresult}, we have $\iota(G^0_X,F) \leq \frac{|E(G^0_X)|+1}{k+2} < \frac{|E(G^0)|+1}{k+2} = 1 + \iota(H_0,F)$, so $G^0_X$ has an $F$-isolating set $D^0_X$ with $|D^0_X| \leq \iota(H_0,F)$. Since $X = N[D_X] \cap V(G^0)$, $D^0_X \cup D_X$ is an $F$-isolating set of $G$, but this contradicts (\ref{isol-components}). Therefore, $G^0_X \simeq F$. Since $E(H_0) \subseteq E(G^0_X)$ and $k = |E(G^0_X)| \geq |E(H_0)| = (k+2)\iota(H_0,F)$, we obtain $|E(H_0)| = 0$, so $V(H_0) = \{y_0\}$ and $E(G^0_X) = (E(F_1) \setminus \{vx\}) \cup \{x_0y_0\}$ with $\{x\} = X$. Thus, $E(G^0) = E(G^0_X) \cup \{vx\}$, and hence $G^0$ is a pure $(k+1,F)$-special graph. By Lemma~\ref{2PureSpecialMulti}, $G$ is a pure $(m,F)$-special graph.   

Now suppose $\mathcal{H}=\{H_0\}$. Then, $G = G^0$. 
Let $\mathcal{J}$ be the set of components of $H_0 - y_0$. Since $H_0$ is connected, for each $J \in \mathcal{J}$, $y_0z_J \in E(G)$ for some $z_J \in V(J)$. Let $\mathcal{J}' = \{J \in \mathcal{J} \colon (J,F) \mbox{ is special}\}$ and $\mathcal{J}^* = \mathcal{J} \setminus \mathcal{J}'$. For each $J \in \mathcal{J}^*$, let $D_J$ be an $F$-isolating set of $J$ of size $\iota(J,F)$. Then, $\{y_0\} \cup \bigcup_{J \in \mathcal{J}^*} D_J$ is an $F$-isolating of $G$. By the induction hypothesis, $|D_J| \leq \frac{|E(J) \cup \{y_0z_J\}|}{k+2}$ for each $J \in \mathcal{J}^*$. Let $h = \sum_{J \in \mathcal{J}^*} |E(J) \cup \{y_0z_J\}|$. We have 
%
%
%
\begin{equation} \frac{m+1}{k+2} = \iota(G,F) \leq 1 + \sum_{J \in \mathcal{J}^*} |D_J| \leq \frac{k + 2 + h}{k+2} = \frac{|E(F_1) \cup \{x_0y_0\}| + 1 + h}{k+2} \leq \frac{m+1}{k+2}, \nonumber
\end{equation}
so $\mathcal{J}' = \emptyset$, $E(G) = E(F_1) \cup \{x_0y_0\} \cup \bigcup_{J \in \mathcal{J}} (E(J) \cup \{y_0z_J\})$, and $|D_J| = \frac{|E(J) \cup \{y_0z_J\}|}{k+2}$ for each $J \in \mathcal{J}$. By the induction hypothesis, for each $J \in \mathcal{J}$, $J$ is a pure $(|E(J)|,F)$-special graph or a copy of $F+e$ for some $e \in E(\Fbar)$. By Lemma~\ref{2PureSpecialMulti}, $G$ is a pure $(m,F)$-special graph. 
\medskip
\\
\textbf{Case 2:} $\mathcal{H}' \neq \emptyset$. For each $x \in N(v)$, let $\mathcal{H}'_x = \{H \in \mathcal{H}' \colon H \mbox{ is linked to } x\}$ and $\mathcal{H}_x^* = \{H \in \mathcal{H} \setminus \mathcal{H}' \colon H \mbox{ is linked to $x$ only}\}$. For each $H \in \mathcal{H} \setminus \mathcal{H}'$, let $D_H$ be an $F$-isolating set of $H$ of size $\iota(H,F)$.

Suppose $|\mathcal{H}'_x| \geq 2$ for some $x \in N(v)$. Let $X = \{x_H \colon H \in \mathcal{H}' \setminus \mathcal{H}'_x\}$. We have $x \notin X \subset N(v)$, so $d(v) \geq 1 + |X|$. Let $D = \{v, x\} \cup X \cup \bigcup_{H \in \mathcal{H} \setminus \mathcal{H}'} D_{H} $. Then, $D$ is an $F$-isolating set of $G$. Since $\iota(G,F) \leq |D|$ and
\begin{align} m + 1 &\geq 1 + d(v) + \sum_{H \in \mathcal{H}'_x} |E(H) \cup \{x y_{x,H}\}| + \sum_{H \in \mathcal{H} \setminus \mathcal{H}_x'} |E(H) \cup \{x_H y_{H}\}| \nonumber \\
&\geq 2 + |X| + (k+1)|\mathcal{H}_x'| + (k+1)|\mathcal{H}' \setminus \mathcal{H}_x'| + \sum_{H \in \mathcal{H} \setminus \mathcal{H}'} (k+2)|D_H| \nonumber \\
&\geq 2 + |X| + 2(k+1) + (k+1)|X| + \sum_{H \in \mathcal{H} \setminus \mathcal{H}'} (k+2)|D_H| = (k+2)|D|, \nonumber
\end{align}
we have $\frac{m+1}{k+2} \geq |D| \geq \iota(G,F) = \frac{m+1}{k+2}$, and hence equality holds throughout. In particular, $|D| = \iota(G,F)$ and $d(v) = 1 + |X|$. The latter equation yields $N(v) = \{x\} \cup X$. Thus, $N(v) \subseteq D$, and hence $D \setminus \{v\}$ is also an $F$-isolating set of $G$. This contradicts $|D| = \iota(G,F)$. 

Therefore, 
\begin{equation} |\mathcal{H}'_x| \leq 1 \mbox{ \emph{for each} } x \in N(v). \label{k=3Hx<2} 
\end{equation} 
Let $H \in \mathcal{H}'$. Let $x = x_H$ and $y = y_{H}$.\medskip
\\
\textbf{Case 2.1:} \emph{$H$ is linked to $x$ only.} Let $X = \{x\} \cup V(H)$. Then, $G - X$ has a component $G_v^*$ such that $N[v] \setminus \{x\} \subseteq V(G_v^*)$, and the other components of $G - X$ are the members of $\mathcal{H}_{x}^*$. Let $D^*$ be an $F$-isolating set of $G_v^*$ of size $\iota(G_v^*,F)$, and let $D = \{x\} \cup D^* \cup \bigcup_{I \in \mathcal{H}_x^*} D_I$. Then, $D$ is an $F$-isolating set of $G$. Since 
\begin{equation} E(G) \supseteq \{vx, xy\} \cup E(H) \cup E(G_v^*) \cup \bigcup_{I \in \mathcal{H}_x^*} (E(I) \cup \{x_Iy_I\}), \label{starpoint1} 
\end{equation}
\begin{equation} m+1 \geq 3 + k + |E(G_v^*)| + \sum_{I \in \mathcal{H}_x^*} (|E(I)| + 1) \geq 3 + k + |E(G_v^*)| + \sum_{I \in \mathcal{H}_x^*} (k+2)|D_I|. \label{starpoint2}
\end{equation}

Suppose $|V(G_v^*)| \leq \ell$. Then, $\{x\} \cup \bigcup_{I \in \mathcal{H}_x^*} D_I$ is an $F$-isolating set of $G$, and hence $\iota(G, F) \leq 1 + \sum_{I \in \mathcal{H}_x^*} |D_I|$. By (\ref{starpoint2}),
\[m + 1 \geq 3 + k + \sum_{I \in \mathcal{H}_x^*} (k+2)|D_I| > (k+2)\iota(G, F),\] 
which contradicts $\iota(G, F) = \frac{m+1}{k+2}$.  

Therefore, $|V(G_v^*)| > \ell$, and hence $F \not\simeq G_v^* \not\simeq F+e$ for each $e \in E(\Fbar)$. By Theorem~\ref{mainresult}, $|D^*| \leq \frac{|E(G_v^*)|+1}{k+2}$, so by (\ref{starpoint2}), 
\begin{equation}
m+1 \geq 3 + k + (k+2)|D^*| - 1 + \sum_{I \in \mathcal{H}_x^*} (k+2)|D_I| = (k+2)|D| \geq (k+2)\iota(G,F). \label{starpoint3}
\end{equation}
Since $\iota(G,F) = \frac{m+1}{k+2}$, we have equality throughout (\ref{starpoint1}), (\ref{starpoint2}) and (\ref{starpoint3}). In particular, $|D^*| = \frac{|E(G_v^*)|+1}{k+2}$ and $|D_I| = \frac{|E(I)|+1}{k+2}$ for each $I \in \mathcal{H}_x^*$. By the induction hypothesis, for each $I \in \{G_v^*\} \cup \mathcal{H}_x^*$, $I$ is a pure $(|E(I)|,F)$-special graph or a copy of $F+e$ for some $e \in \Fbar$. Moreover, $G[X]$ is a $(k+1,F)$-special graph. By Lemma~\ref{2PureSpecialMulti}, $G$ is a pure $(m,F)$-special graph.\medskip
\\
\textbf{Case 2.2:} \emph{$H$ is linked to some $x' \in N(v) \setminus \{x\}$.} Then, $x'y' \in E(G)$ for some $y' \in V(H)$. Let $J = G - V(H)$. Then, $J$ is connected. Since $H \simeq F$, $V(H) \subseteq N[w]$ for some $w \in V(H)$. Let $A = E_G(N(v), V(H))$. Then, $xy, x'y' \in A$.

Suppose that $(J, F)$ is special. Then, $J \simeq F$, and hence $J = F_1$. Since $V(F_1) \subseteq N[v]$ and $V(H) \subseteq N[w]$, $V(G) = N[\{v, w\}]$. We have $\frac{m+1}{k+2} = \iota(G, F) \leq 2$, contradicting $m \geq 2k+4$.

Therefore, $(J, F)$ is not special. By Lemma~\ref{lemma} (with $X = \{w\}$ and $Y = V(H)$) and the induction hypothesis, 
\[\iota(G,F) \leq 1 + \iota(J,F) \leq \frac{|E(H) \cup A|}{k+2} + \frac{|E(J)|+1}{k+2} \leq \frac{m+1}{k+2}.\]
Since $\iota(G,F) = \frac{m+1}{k+2}$, we have $\iota(G,F) = 1 + \iota(J,F)$, $A = \{xy,x'y'\}$, $\iota(J,F) = \frac{|E(J)|+1}{k+2}$ and $m = |E(H) \cup A| + |E(J)| = k + 2 + |E(J)|$. Since $m \geq 2k+4$, $|E(J)| \geq k+2$. By the induction hypothesis, $J$ is an $(|E(J)|,F)$-special graph.

Suppose that $x$ and $x'$ are in distinct $F$-constituents of $J$. By Lemma~\ref{MinIsolatingSets2vertices}, $J$ has a minimum $F$-isolating set $D_J$ such that $x, x' \in D_J$. Thus, $D_J$ is an $F$-isolating set of $G$. This contradicts $\iota(G,F) = 1 + \iota(J,F)$.

Therefore, $x$ and $x'$ are in the same $F$-constituent of $J$. Let $J$ have quotient graph $T$, $F$-connections $v_1, \ldots, v_q$, and $F$-constituents $G_1, \ldots, G_q$, where $v_i \in V(G_i)$ for each $i \in [q]$. We may assume that $x, x' \in V(G_q)$. Since $|E(J)| \geq k+2$, $q \geq 2$. Let $D_J = \{v_i \colon i \in [q-1]\}$. Since $T$ is connected, $v_q \in N[D_J]$. Let $I = G[V(G_q - v_q) \cup V(H)]$. Since $x \neq x'$, we may assume that $x \neq v_q$. Thus, $xy \in E(I)$. We have $H \simeq F \simeq G_q - v_q$, so $I$ is connected and $|E(I)| \leq 2k + |A| = 2k+2$. By Theorem~\ref{mainresult}, $I - N_I[x_I]$ is $F$-free for some $x_I \in V(I)$. Since $v_q \in N[D_J]$, $\{x_I\} \cup D_J$ is an $F$-isolating of $G$. Thus, we have $\iota(G,F) \leq q = \frac{|E(J)|+1}{k+2} < \frac{m+1}{k+2}$, contradicting $\iota(G,F) = \frac{m+1}{k+2}$.~\hfill{$\Box$}
\\
\\
\textbf{Acknowledgements:} The authors are grateful to an anonymous referee for checking the paper carefully and providing constructive remarks.

\end{document}